\documentclass[12pt,twoside]{article}
\usepackage{amsmath}
\usepackage{amsfonts}
\usepackage{amsthm}
\usepackage{amssymb}
\usepackage{graphicx}
\usepackage{multicol}

\begin{document}
\title{{\normalsize{\bf Another proof of free ribbon lemma}}}
\author{{\footnotesize Akio Kawauchi}\\
{\footnotesize{\it Osaka Central Advanced Mathematical Institute, Osaka Metropolitan University}}\\
{\footnotesize{\it Sugimoto, Sumiyoshi-ku, Osaka 558-8585, Japan}}\\
{\footnotesize{\it kawauchi@omu.ac.jp}}}
\date\,
\maketitle
\vspace{0.25in}
\baselineskip=10pt
\newtheorem{Theorem}{Theorem}[section]
\newtheorem{Conjecture}[Theorem]{Conjecture}
\newtheorem{Lemma}[Theorem]{Lemma}
\newtheorem{Sublemma}[Theorem]{Sublemma}
\newtheorem{Proposition}[Theorem]{Proposition}
\newtheorem{Corollary}[Theorem]{Corollary}
\newtheorem{Claim}[Theorem]{Claim}
\newtheorem{Definition}[Theorem]{Definition}
\newtheorem{Example}[Theorem]{Example}

\begin{abstract}  Free ribbon lemma  that every free sphere-link in the 4-sphere is a ribbon sphere-link is shown in an earlier paper by the author. In this paper, another proof of this lemma  
is given. 

\phantom{x}



\end{abstract}

\baselineskip=15pt

\bigskip

\noindent{\bf 1. Introduction}

A {\it surface-link} is a closed oriented (possibly, disconnected) surface $L$ smoothly embedded in the 
4-sphere $S^4$. When $L$ is connected, $L$ is called a {\it surface-knot}. 
If $L$ consists of 2-spheres $L_i\,(i=1,2,\dots,n)$, then $L$ is called 
a {\it sphere-link} (or an $S^2$-{\it link}) of $n$  components. 
It is shown that a surface-link $L$  is a trivial surface-link (i.e., bounds disjoint handlebodies in $S^4$) if $\pi_1(S^4\setminus L,x_0)$ is a meridian-based free group (see  \cite{K21,K23-1, K23-2}).  
A surface-link $L$ is {\it ribbon} if $L$ is obtained from a trivial $S^2$-link $O$ in $S^4$ 
by surgery along smoothly embedded disjoint 1-handles on $O$. 
A surface-link $L$ in the 4-sphere $S^4$ is {\it free} if the fundamental group 
$\pi_1(S^4\setminus L,x_0)$ is a (not necessarily meridian-based) free group.  
The {\it free ribbon lemma} is the following theorem.

\phantom{x}

\noindent{\bf Theorem.} Every free $S^2$-link is a ribbon $S^2$-link. 

\phantom{x}

This theorem is a basic result concerning Whitehead aspherical conjecture \cite{K24-1} \cite{K24-3} and classical Poincar{\'e} conjecture \cite{K24-2}, and the proof is done in \cite{K24-1} as an appendix.  {\it At present, it appears unknown whether or not every free surface-link is a ribbon surface-link. }
In this paper, another proof of this theorem is given as follows. 

\phantom{x}

\noindent{\bf Proof of Theorem.} 
Let $L_i\,(i=1,2,\dots,n)$ be the components of a free $S^2$-link $L$. 
Let $x_i,\,(i=1,2,\dots,n)$ be a basis of the free fundamental group 
$G=\pi_1(S^4\setminus L,x_0)$. 
Let $y_i$ be a meridian element of $L_i$ in $G$, so that $y_i\,(i=1,2,\dots,n)$ are a meridian 
system of $G$.
By Nielsen transformations, $y_i$ is equal to $x_i$ modulo the commutator subgroup $[G,G]$ of 
$G$. It is known that the group $G$ is isomorphic to a group $G^P$ with  
 Wirtinger presentation 
\[P=<y_{ij}\, (1\leq j\leq m_i,\, 1\leq i\leq  n )|\, r_{ij}\, (2\leq j\leq m_i+s_i,\, 1\leq i\leq n)>\]
such that $y_{i1}=y_i\, (i=1,2,\dots,n)$ and the relators 
$r_{ij}\, (j=2,3,\dots, m_i+s_i,\, i=1,2,\dots, n)$ 
are given by 
$r_{ij}:\, y_{ij}=w_{ij}y_{i 1}w_{ij}^{-1} $ for $j$ with $2\leq j\leq m_i$, $1\leq i\leq n$, and 
$r_{ij}:\, y_{i1}=w_{ij}y_{i 1}w_{ij}^{-1}$ for $j$ with $m_i+1\leq j\leq m_i+s_i$, $1\leq i\leq n$, 
where $w_{ij}\, (j=2,3,\dots, m_i+s_i,\, i=1,2,\dots, n)$ are words in the letters 
$y_{ij}\, (j=1,2,\dots, m_i,\, i=1,2,\dots, n )$. This result is obtained from Yajima \cite{Yajima70} 
because $G$ has a weight system $y_i\,(i=1,2,\dots,n)$, $ H_1(G;Z)\cong Z^n$ and $H_2(G;Z)=0$. 
It is observed that this result can be also obtained by an alternative geometric method using 
a normal form of a surface-link in $R^4$ \cite{KSS}. 
In fact, put the $S^2$-link $L$ in a normal form of in the 4-space $R^4$ with 
$L[0]=L\cap R^3[0]$ a middle cross-sectional link and calculate the fundamental groups 
$\pi_1(R^3[0,+\infty)\setminus L\cap R^3[0,+\infty),x_0)$ and 
$\pi_1(R^3(-\infty,0]\setminus L\cap R^3(-\infty,0],x_0)$ with Wirtinger presentations starting from the fundamental group $\pi_1(R^3[0]\setminus L[0],x_0)$ with a Wirtinger presentation 
to obtain the  group $G$ with a Wirtinger presentation by van Kampen theorem. 
See \cite{K07, K08} for this construction and \cite{Kamada01} for a generalization. 
By fixing an isomorphism $G^P\to G$, regard the generators $y_{ij}\, (j=1,2,\dots, m_i,\, i=1,2,\dots, n )$ of $P$ as fixed words in the basis $x_i\,(i=1,2,\dots,n)$ of $G$. 
Then the relator $y_{i1}=w_{ij}y_{i1}w_{ij}^{-1}$ for every $i$ and  $j$ with 
$m_i+1\leq j\leq m_i+s_i$ can be written as $y_{i1}=a_{ij}^{u(i,j)}$ and $w_{ij}=a_{ij}^{v(i,j)}$ for a word $a_{ij}$ in $x_i\,(i=1,2,\dots,n)$ and some integers 
$u(i,j), v(i,j)$ because any nontrivial abelian subgroup of a free group is an infinite cyclic group. 
The elements $y_i=y_{i1}\, (i=1,2,\dots,n)$ form the same abelian basis as  $x_i\, (i=1,2,\dots,n)$ in the free abelian group $G/[G,G]$, so that $u(i,j)=\pm1$ for every $i$ and  $j$. 
Thus, $w_{ij}=y_{i1}^{u(i,j)v(i,j)}$ for every $i$ and  $j$ with $m_i+1\leq j\leq m_i+s_i$, which means 
that the relators 
$r_{ij:}\, y_{i1}=w_{ij}y_{i1}w_{ij}^{-1}\, (m_i+1\leq j\leq m_i+s_i)$ 
are identity relations in the free group $<y_{ij}\, (1\leq j\leq m_i,\, 1\leq i\leq  n  )>$.
Thus, the Wirtinger presentation $P$ is equivalent to the Wirtinger presentation 
\[R=<y_{ij}\, (1\leq j\leq m_i,\, 1\leq i\leq  n  )|\, r_{ij}\, (2\leq j\leq m_i,\, 1\leq i\leq  n )>\]
with $y_{i1} =y_i \, (i=1,2,\dots, n)$ and the relators 
$r_{ij}\, (2\leq j\leq m_i,\, 1\leq i\leq  n )$ given by 
$r_{ij}:\, y_{ij}=w_{ij}y_{i1}w_{ij}^{-1} \,(2\leq j\leq m_i,\, 1\leq i\leq  n )$.
By Yajima's construction \cite{Yajima62} (see also \cite{K07,K08}), 
there is a ribbon $S^2$-link $L^R$ with the fundamental group 
$G^R=\pi_1(S^4\setminus L^R,x_0)$ of the Wirtinger presentation $R$
which is isomorphic to $G$ by an isomorphism $G^R\to G$  
sending  a meridian element $y^R_i$ of the $i$th component $L^R_i$ of $L^R$ to the meridian element $y_i$ of $L_i$ in $G$ for every $i\,(i=1,2,\dots,n)$ and a basis $x^R_i\,(i=1,2,\dots,n)$ of $G^R$ to the basis $x_i\,(i=1,2,\dots,n)$ of $G$. 
Let $Y^R$ and $Y$ be the 4D manifolds (both diffeomorphic to the $n$-fold connected sum 
of $S^1\times S^3$) obtained from $S^4$ by surgeries along $L^R$ and $L$, 
respectively, and $\ell^R_i\,(i=1,2,\dots,n)$ and $\ell_i\,(i=1,2,\dots,n)$ the loop systems 
obtained from $L^R_i\,(i=1,2,\dots,n)$ and $L_i\,(i=1,2,\dots,n)$, respectively. 
By \cite{K24-1}, there is an orientation-preserving diffeomorphism $f:Y^R\to Y$ 
sending the loop system $\ell^R_i \,(i=1,2,\dots,n)$ to the loop system $\ell_i \,(i=1,2,\dots,n)$. 
Note that this result is obtained from the  smooth unknotting conjecture for $S^2$-knots \cite{K21,K23-1, K23-2} and the 4D smooth Poincar{\'e} conjecture \cite{K23-3}. 
By the back surgeries from $Y^R$ to $S^4$ 
along $\ell^R_i \,(i=1,2,\dots,n)$ and from $Y$ to $S^4$ along $\ell_i \,(i=1,2,\dots,n)$, this diffeomorphism $f$ induces an orientation-preserving diffeomorphism $f':S^4\to S^4$ sending $L^R$ to $L$. Thus, the $S^2$-link $L$ is a ribbon $S^2$-link. 
This completes the proof of Theorem. 

\phantom{x}

In the proof of Theorem, the ribbon $S^2$-link $L^R$ is called a {\it ribbon presentation} 
of the free $S^2$-link $L$. 
The following corollary is obtained from the proof of Theorem. 

\phantom{x}

\noindent{\bf Corollary.} Let $L$ be a free $S^2$-link 
in the 4-sphere $S^4$ containing a free $S^2$-link $K$ as a sublink. 
For any ribbon presentation of $K^R$ of $K$, there is a ribbon presentation $L^R$ of $L$ 
containing $K^R$ as a sublink. 

\phantom{x}

\noindent{\bf Proof of Corollary.} 
The ribbon presentation of $K^R$ of $K$ is in a normal form. Thus, the result is obtained 
from the observation that a normal form of $L$ is taken to contain $K^R$ as a sublink (see \cite{KSS}). 

\phantom{x}

\noindent{\bf Acknowledgements.}  
This work was partly supported by JSPS KAKENHI Grant Number JP21H00978 and MEXT Promotion of Distinctive Joint Research Center Program JPMXP0723833165. 

\phantom{x}

\end{document}